\documentclass[preprint,12pt]{elsarticle}
\usepackage{setspace}
\usepackage[usenames]{color}
\usepackage{tikz}
\usetikzlibrary{calc,quotes,angles,positioning, arrows}
\usepackage{graphicx}
\usepackage{soul}
\usepackage{abraces}
\usepackage{multicol}
\usepackage{enumerate}
\usepackage{stackengine}

\usepackage[font=footnotesize,skip=0pt]{caption}

\usepackage[colorlinks=true,
linkcolor=webgreen,
filecolor=webbrown,
citecolor=webgreen]{hyperref}

\definecolor{webgreen}{rgb}{0,.5,0}
\definecolor{webbrown}{rgb}{.6,0,0}
\definecolor{limegreen}{rgb}{0.2, 0.8, 0.2}

\usepackage{tcolorbox}
\usepackage{color}
\usepackage{fullpage}
\usepackage{float}

\usepackage{graphics}
\usepackage{amsthm}
\usepackage{amsmath}
\usepackage{amsfonts}
\usepackage{latexsym}
\usepackage{epsf}
\usepackage{framed,enumitem}

\journal{Australasian Journal of Combinatorics}

\begin{document}

\newtheorem{theorem}{Theorem}[section]
\newtheorem{lemma}[theorem]{Lemma}
\newtheorem{corollary}[theorem]{Corollary}
\newtheorem{proposition}[theorem]{Proposition}

\theoremstyle{definition}
\newtheorem{defn}[theorem]{Definition}
\newtheorem{exmp}[theorem]{Example}
\newtheorem{remark}[theorem]{Remark}
\newtheorem{algo}[theorem]{Algorithm}

\begin{frontmatter}

\title{Permutation-generated maps between Dyck paths}

\author[inst1]{Kevin Limanta\corref{cor1}}
\ead{k.limanta@unsw.edu.au}
\cortext[cor1]{Corresponding author.}

\affiliation[inst1]{organization={School of Mathematics and Statistics},
            addressline={University of New South Wales Sydney},
            country={Australia}}

\author[inst2]{Hopein Christofen Tang}
\ead{hopeinct@students.itb.ac.id}

\affiliation[inst2]{organization={Department of Mathematics},
            addressline={Institut Teknologi Bandung}, 
            country={Indonesia}}
            
\author[inst3]{Yozef Tjandra}
\ead{yozef.tjandra@calvin.ac.id}

\affiliation[inst3]{organization={Department of IT and Big Data Analytics},
            addressline={Calvin Institute of Technology},
            country={Indonesia}}
        
\begin{abstract}
In 2003, Deutsch and Elizalde defined a family of bijective maps between the set of Dyck paths to itself which is induced by some particular permutations. In this paper, we extend the construction of the maps by allowing the permutation to be arbitrary. We characterise the permutations which generate the same map and find all permutations generating a bijection among Dyck paths. Consequently, we give a new combinatorial interpretation of the quantity $(2n-1)!!$ as well as some new statistics of Dyck paths which are equidistributed to some known height statistics via our generalised maps.
\end{abstract}


\end{frontmatter}

\section{Introduction}
Bijective methods are one of the classic tools in enumerative combinatorics. One enumeration problem can be translated to another, often in a different discrete structure.

There are various results about bijections between the set of Dyck paths to itself as well as to some other discrete structures. For example, Deutsch showed that the number of high peaks has the Narayana distribution  \cite{deutsch1998bijection}. He also proved in \cite{deutsch1999involution} that the number of valleys and the number of double-rises, as well as the height of the first peak and the number of returns are equidistributed. Furthermore, there are studies about bijections from the set of Dyck paths to pattern-avoiding permutations \cite{elizalde2011fixed,elizalde2003simple,elizalde2004bijections} and to bar graphs \cite{deutsch2018bijection}. 

In 2003, Deutsch and Elizalde defined a family of bijections between the set of Dyck paths to itself in \cite{elizalde2003simple}. The maps, deemed by the authors as simple and unusual, are generated by some particular permutation in the symmetric group and have a pleasant property that they establish an equidistribution between one statistic of Dyck paths to other known statistics as well as some results in enumerating pattern-avoiding permutations. 

In this paper, we work on a natural generalisation of the maps defined by Deutsch and Elizalde by letting the maps be generated by arbitrary permutations. We first noticed that some maps produced by different permutations might be completely identical. From this, we studied the complete characterisation of such permutations. Furthermore, among those permutation-generated maps, we find a much larger family of bijections compared to that in \cite{elizalde2003simple}. We call any permutation generating such bijection a \emph{circularly-connected permutation} (CCP).

Using our generalised maps, we provide an alternative combinatorial interpretation of an identity involving $(2n-1)!!$ from Callan \cite{callan2009combinatorial}. Moreover, we obtain a connection between a Dyck path enumeration problem with respect to some existing height statistics and the number of \emph{unmatched} steps. The latter is a new statistic whose meaning will be introduced in Section 5. We show that the distribution of those statistics among all Dyck Paths are equal by utilising our bijections.

This paper is organised as follows. The notation and basic definitions are provided in Section 2. In Section 3, we define the generalised maps and their intermediate consequences. In Section 4, we prove that CCPs and only CCPs generate bijections. In the last section, we discuss some results on statistics on Dyck paths.

\section{Dyck paths: generalities and terminologies}

We adopt the following notation of which we also summarise in Appendix B. We denote the set of integers from $a$ to $b$ inclusive by $[a,b]$. For brevity, $[n] := [1,n]$. We  always assume that any arithmetic operation imposed on $[n]$ to be modular. The group of permutations of elements of $[n]$ is denoted by $S_n$. Following the notation in \cite{elizalde2003simple}, for $\sigma \in S_{n}$ we write $\sigma_k$ in place of $\sigma(k)$. By $\sigma_{[k]}$, we mean $\{\sigma_1, \sigma_2, \ldots, \sigma_k\}$.

Given a two-dimensional integer lattice $\mathbb{Z}^2$, we can make a lattice path starting at $(0,0)$ consisting of up-steps $(1,1)$ and down-steps $(1,-1)$, represented by $\mathtt{u}$ and $\mathtt{d}$ respectively. A \emph{Dyck path} of size $n$ is such a lattice path between $(0,0)$ and $(2n,0)$ which never goes below the $x$-axis. Any Dyck path can be encoded by a string of $\mathtt{u}$ and $\mathtt{d}$ which we call its \emph{Dyck word}. Step $k$ of a Dyck path $P$ is denoted by $P_k \in \{\mathtt{u}, \mathtt{d}\}$. We will use the notation $\mathtt{u}^k$ to indicate a string of $k$ consecutive up-steps with $\mathtt{d}^k$ defined similarly. Finally, $\mathcal{D}_n$ denotes the set of all Dyck paths of size $n$.  

\begin{figure}[H]
\begin{center}
\begin{tikzpicture}[scale=0.5]
\tikzstyle{every node}=[draw,circle,fill=black,minimum size=2pt,inner sep=0pt,radius=1pt]
      \node (A) at (0,0) {};
      \node (B) at (1,1) {};
      \node (C) at (2,2) {};
      \node (D) at (3,1) {};
      \node (E) at (4,2) {};
      \node (F) at (5,3) {};
      \node (G) at (6,2) {};
      \node (H) at (7,1) {};
      \node (I) at (8,0) {};
      \draw (A)--(B);
      \draw (B)--(C);
      \draw (C)--(D);
      \draw (D)--(E);
      \draw (E)--(F);
      \draw (F)--(G);
      \draw (G)--(H);
      \draw (H)--(I);
\end{tikzpicture}
\end{center}
\caption{A Dyck path $P\in\mathcal{D}_4$ with Dyck word $\mathtt{uuduuddd}=\mathtt{u}^2\mathtt{du}^2\mathtt{d}^3$, $P_2 = \mathtt{u}$, and $P_7 = \mathtt{d}$.}
\label{Fig1}
\end{figure}

\vspace{-4mm}
For any $D \in \mathcal{D}_n$, following Elizalde's original definition in \cite{elizalde2003simple}, we define a \emph{tunnel} of $D$ to be the horizontal segment between two lattice points of $D$ that intersects $D$ only in those two points and always stays below $D$. If $D \in \mathcal{D}_n$, then $D$ has exactly $n$ tunnels and each tunnel can be associated with a pair of $(k,l)$ such that step $k$ and step $l$ are connected by a tunnel.
\begin{figure}[H]
\begin{center}
\begin{tikzpicture}[scale=0.5]
\tikzstyle{every node}=[draw,circle,fill=black,minimum size=2pt,inner sep=0pt,radius=1pt]
      \node (A) at (0,0) {};
      \node (B) at (1,1) {};
      \node (C) at (2,2) {};
      \node (D) at (3,1) {};
      \node (E) at (4,2) {};
      \node (F) at (5,3) {};
      \node (G) at (6,2) {};
      \node (H) at (7,1) {};
      \node (I) at (8,0) {};
      \draw (A)--(B);
      \draw (B)--(C);
      \draw (C)--(D);
      \draw (D)--(E);
      \draw (E)--(F);
      \draw (F)--(G);
      \draw (G)--(H);
      \draw (H)--(I);
      \draw[red,thin,dashed] (A) -- (I);
      \draw[blue,thin,dashed] (B) -- (D);
      \draw[red,thick,dashed] (D) -- (H);
      \draw[blue,thick,dashed] (E) -- (G);
       \end{tikzpicture}
\end{center}
\caption{The four tunnels of $\mathtt{uuduuddd}$.}
\label{Fig2}
\end{figure}

\vspace{-4mm}
The tunnel pairs of $D$ can be encoded as a permutation $\tau_D$, called the \emph{tunneling} of $D$, where $\tau_D(k) = \ell$ if and only if $(k,\ell)$ is a tunnel pair. For example, if $D$ is the Dyck path in Figure \ref{Fig2} above, then the four tunnels are the horizontal line segments connecting $(0,0)$ and $(8,0)$, $(1,1)$ and $(3,1)$, $(3,1)$ and $(7,1)$, as well as $(4,2)$ and $(6,2)$ with $\tau_D = 83276541$.

A Dyck path $D \in \mathcal{D}_n$ can also be represented as a circle labeled with points in $[2n]$ arranged in a clockwise manner, where vertices with label $k$ and $\ell$ are joined by a chord if and only if $\tau_D(k) = \ell$. We call this the \emph{circular representation} of the Dyck path $D$. For example, the above Dyck path admits the following circular representation. 

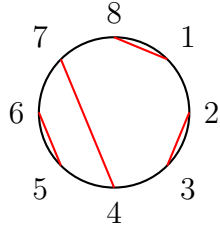
\vspace{-2mm}

\begin{figure}[H]
    \centering
    \begin{tikzpicture}[scale=0.5]
        \draw[thick] (0,0) circle (2cm);
        \coordinate (A) at (45:2);
        \node at (A) [above right = 0.3535mm of A] {$1$};
        \coordinate (B) at (0:2);
        \node at (B) [right = 0.5mm of B] {$2$};
        \coordinate (C) at (315:2);
        \node at (C) [below right = 0.3535mm of C] {$3$};
        \coordinate (D) at (270:2);
        \node at (D) [below = 0.5mm of D] {$4$};
        \coordinate (E) at (225:2);
        \node at (E) [below left = 0.3535mm of E] {$5$};
        \coordinate (F) at (180:2);
        \node at (F) [left = 0.5mm of F] {$6$};
        \coordinate (G) at (135:2);
        \node at (G) [above left = 0.3535mm of G] {$7$};
        \coordinate (H) at (90:2);
        \node at (H) [above = 0.5mm of H] {$8$};
        \draw[red,thick] (A) -- (H);
        \draw[red,thick] (B) -- (C);
        \draw[red,thick] (D) -- (G);
        \draw[red,thick] (E) -- (F);
    \end{tikzpicture}
\end{figure}
\vspace{-4mm}
\noindent It is clear that the $n$ chords are non-intersecting since $D$ is a Dyck path.

Given a Dyck path $P$, let $h(P)$ be the height of the highest peak in $P$ and $h_k(P)$ denote the height of $P$ after step $k$, that is $k$ subtracted by twice the number of down-steps until step $k$. For example, if $P$ is the Dyck path in Figure \ref{Fig1}, then $h_2(P)=2, h_3(P)=1, h_5(P)=3, h_7(P)=1$.

\section{Construction of the map and some consequences}

\subsection{\texorpdfstring{$\sigma$-paths}{sigma-paths}}
We are ready to define our permutation-generated map.
\begin{defn}\label{Def1}
    Given a permutation $\sigma \in S_{2n}$ and a Dyck path $D$, the $\sigma$-path of $D$, denoted by $\sigma(D)$, is a lattice path constructed by the following algorithm: at iteration $k \in [2n]$, 
    \[
    \sigma(D)_k =\begin{cases}
        \texttt{u} & \text{ if } \tau_D(\sigma_k) \notin \sigma_{[k]},\\
        \texttt{d} & \text{ otherwise.}
    \end{cases}
    \]
    The map $D \mapsto \sigma(D)$ is denoted by $\sigma(\cdot)$.
\end{defn}
\begin{remark}
    The condition $\tau_D(\sigma_k) \notin \sigma_{[k]}$ is equivalent to $\sigma^{-1}\tau_D(\sigma_k)>k$.
\end{remark}

\begin{remark}
A more informal way to describe the construction of $\sigma(D)$ is as follows: at iteration $k \in [2n]$, step $k$ of $\sigma(D)$ is \texttt{u} if and only if the tunnel pair of step $\sigma_k$ of $D$ has not been read before. 
\end{remark}
The readers are encouraged to compare this with the construction of the map $\Phi$ by Deutsch and Elizalde in \cite[Section 3]{elizalde2003simple}. Indeed, the map $\Phi$ in their construction is  $\sigma(\cdot)$ with \[\sigma_k=
\begin{cases}
\frac{k+1}{2} &\mbox{if $k$ is odd,}\\
2n+1-\frac{k}{2} &\mbox{if $k$ is even.}
\end{cases}\]

It is not hard to see that $\sigma(D)$ is a Dyck path for all $\sigma \in S_{2n}$ and $D \in \mathcal{D}_n$. The resulting word of $\sigma(D)$ has precisely $n$ up-steps and down-steps since $D$ has $n$ tunnels. Moreover, the number of down-steps in the first $k$ positions will not be less than that of up-steps for any $k\in[2n]$ since the only way a down-step in position $k$ is produced is to have the tunnel pair of step $\sigma_k$ of $D$ read. Before it is read, an up-step must have been written at some step $\ell < k$. Hence the resulting word is a Dyck word.

\begin{exmp}
Let $D \in \mathcal{D}_4$ be represented by its Dyck word $\mathtt{uuddudud}$ and $\sigma = 14285763$. Now,
\begin{itemize}
    \item at $k = 1$, $\tau_D(\sigma_1) = 4 \notin \sigma_{[1]} = \{1\}$, so $\sigma(D)_1 = \mathtt{u}$,
    \item at $k = 2$, $\tau_D(\sigma_2) = 1 \in \sigma_{[2]} = \{1, 4\}$ so $\sigma(D)_2 = \mathtt{d}$,
    \item at $k = 3$, $\tau_D(\sigma_3) = 3 \notin \sigma_{[3]} = \{1, 4, 2\}$ so $\sigma(D)_3 = \mathtt{u}$,
\end{itemize}
and so on. Completing the process, the Dyck word of $\sigma(D)$ is $\mathtt{uduuuddd}$.
\end{exmp}

\subsection{Number of different permutation-generated maps}

A quick observation suggests that two different permutations could correspond to the same map. For instance, consider the permutations $\lambda = 1234$ and $\mu = 2143$. From this, it is natural to investigate the conditions for which two permutations admit the same map as well as the number of different permutation-generated maps. It is easy to check that there is exactly one map when $n = 1$ and three different maps when $n = 2$. The enumeration of such maps for $n \ge 3$ is non-trivial. We present our first result below.

\begin{theorem}\label{ThmCl1}
    For $n\ge 3$, there are
    \begin{align*}
    1-n^2+2\sum_{a=1}^{n-1}\sum_{b=1}^{n-1} \frac{n!n!}{\mathrm{max}\{a,2\}!\mathrm{max}\{b,2\}!}\Bigg({{2n-2-a-b}\choose{n-2}}+{{2n-2-a-b}\choose{n-1-a}}\Bigg)
\end{align*}
different permutation-generated maps.
\end{theorem}

To prove the result above, we need to define some terminologies first. Define an equivalence relation $\sim$ on $S_{2n}$ with $\lambda \sim \mu$ if and only if $\lambda(D) = \mu(D)$ for any $D \in \mathcal{D}_n.$ The equivalence class of $\sigma$ will be denoted by $\mathrm{class}(\sigma)$. It is easy to see that the problem of counting the number of different maps is equivalent to that of counting the number of distinct $\mathrm{class}(\sigma)$. This is done by giving sufficient and necessary conditions for two permutations to be in the same class. 

\begin{defn}
For $n \ge 3$, the \emph{family} of a permutation $\lambda\in S_{2n}$ is the set
\[
\mathrm{fam}(\lambda):=\{\mu \in S_{2n} \,\colon\, \lambda_{[2]}=\mu_{[2]} \text{ and } \lambda_i=\mu_i \text{ for all } i \in [3, 2n-2]\}.\]
\end{defn}

For example, $\mathrm{fam}(12345678)=\{12345678, 12345687, 21345678, 21345687\}$. Note that $\mathrm{fam}(\lambda)$ always has exactly four elements and induces a partition on $S_{2n}$.

\begin{proposition}[Sufficient condition 1]\label{Prop1}
Two permutations that belong to the same family belong to the same class.
\end{proposition}

\begin{proof}
    For any $\mu \in \mathrm{fam}(\lambda)$ and $D \in \mathcal{D}_n$, we have $\mu(D)_2=\mathtt{u}$ if $\lambda(D)_2=\mathtt{u}$ (since $\tau_D(\lambda_1) \neq  \lambda_2$) and $\mu(D)_2=\mathtt{d}$ if $\lambda(D)_2=\mathtt{d}$ (since $\tau_D(\lambda_1)= \lambda_2$). Thus, $\lambda(D)_2=\mu(D)_2$ for all $D \in \mathcal{D}_n.$ Similarly, $\lambda(D)_{2n-1}=\mu(D)_{2n-1}$ for all $D \in \mathcal{D}_n.$ Since $\lambda_i=\mu_i$ for all $i \in [3, 2n-2]$, we have $\lambda(D)=\mu(D)$ for all $D \in  \mathcal{D}_n$. We conclude that $\mathrm{class}(\lambda)=\mathrm{class}(\mu)$ for all $\mu \in \mathrm{fam}(\lambda)$.  
\end{proof}

\begin{defn}
For $n \ge 3$, the \emph{parity} of a permutation $\sigma$ is a pair of positive integers $\mathrm{par}(\sigma) = (a,b)$ such that $a$ is the smallest integer with $\sigma_a\not\equiv \sigma_{a+1}\mod 2 $ and $b$ is the smallest integer with $\sigma_{2n+1-b}\not\equiv \sigma_{2n-b}\mod 2$. 
\end{defn}

For example, $\mathrm{par}(13275468)=(2,3)$. Observe that for any $D \in \mathcal{D}_n$, $\mathrm{par}(\sigma)=(a,b)$ implies that $\sigma(D)_i=\mathtt{u}$ for all $i \in [a]$ and $\sigma(D)_{j}=\mathtt{d}$ for all $j \in [2n+1-b, 2n]$.

\begin{proposition}[Necessary condition 1]\label{Prop2}
Two permutations that belong to the same class have the same parity.
\end{proposition}

\begin{proof}
    Let $\mathrm{par}(\lambda)=(a_1,b_1)$ and $\mathrm{par}(\mu)=(a_2,b_2)$. Suppose that $\mathrm{class}(\lambda)=\mathrm{class}(\mu)$ but $\mathrm{par}(\lambda)\neq\mathrm{par}(\mu)$. Without loss of generality, let $a_1\neq a_2$ with $a_1\geq a_2+1$. Obviously, there exists a Dyck path $D \in \mathcal{D}_n$ with $\tau_D(\mu_{a_2})=\mu_{a_2+1}$ since $\mu_{a_2}\not\equiv \mu_{a_2+1}\,\mathrm{mod}\, 2 $. In this case, $\mu(D)_{a_2+1}=\mathtt{d}$. On the other hand, $a_2+1\leq a_1$ implies $\lambda_i \equiv \lambda_{a_2+1} \,\mathrm{mod}\, 2$ for all $i \in [a_2+1]$. Thus, $\lambda(D)_{a_2+1}=\mathtt{u}$. We have $\lambda(D)\neq \mu(D)$, a contradiction. 
\end{proof}

\begin{proposition}[Necessary condition 2]\label{Prop3}
    Suppose $\mathrm{par}(\lambda) = \mathrm{par}(\mu)$. If there exists three distinct integers $1\leq i<j<k\leq 2n$ and distinct integers $P, Q, R \in [2n]$ such that 
    \begin{enumerate}
        \item $Q \not\equiv R \mod 2$,
        \item $P+1 \neq Q,R$,
        \item $\{\lambda_i,\lambda_j\}=\{P,P+1\}, \{\mu_j,\mu_k\}=\{Q,R\}$ or
        $\{\lambda_i,\lambda_j\}=\{Q,R\},\{\mu_j,\mu_k\}=\{P,P+1\}$,
    \end{enumerate}
    then $\mathrm{class}(\lambda) \neq \mathrm{class}(\mu)$.
\end{proposition}
\begin{proof}
    Suppose such $i,j,k$ exist. Obviously, there exists a Dyck path $D \in \mathcal{D}_n$ such that $\tau_D(P)=P+1$ and $\tau_D(Q)=R$. In other words, $\tau_D(\lambda_i)=\lambda_j$ and $\tau_D(\mu_j)=\mu_k$. Since $i<j<k$, we have $\lambda(D)_j=\mathtt{d}$ and $\mu(D)_j=\mathtt{u}$. This shows that $\lambda(D) \neq \mu(D)$, so $\mathrm{class}(\lambda) \neq \mathrm{class}(\mu)$. 
\end{proof}

To illustrate Proposition \ref{Prop3}, take $\lambda = 13275468$ and $\mu = 53271468$ as examples. The triples $(i,j,k)=(3,5,8)$ and $(P,Q,R) = (8,2,5)$ satisfy all conditions in the proposition.

\begin{defn}
Suppose $\mathrm{par}(\lambda) = \mathrm{par}(\mu) = (a,b)$. We call $\lambda$ and $\mu$ \emph{friends} precisely when the following conditions are satisfied:
\begin{enumerate}
    \item $\lambda_{[a]}=\mu_{[a]}$
    \item $\lambda_i=\mu_i$ for all $i \in [a+1,2n-b]$.
\end{enumerate}
In other words, $\lambda$ and $\mu$ are friends if the first $a$ and last $b$ entries of $\lambda$ and $\mu$ are permutation of each other and their middle entries are identical.
\end{defn}

For example, $13275468$ and $31275846$ are friends while the pair $13275468$ and $53271468$ are not friends. The following proposition is obvious.

\begin{proposition}[Sufficient condition 2]\label{Prop4}
    Two permutations that are friends belong to the same class.
\end{proposition}

Propositions \ref{Prop1} - \ref{Prop4} are the ingredients to prove the following characterisation theorem. The complete proof is by casework and tedious, so it will be presented in Appendix A. A big part of the proof is to establish the existence of a triple $(i,j,k)$ that satisfies Proposition \ref{Prop3} in every possible case to show that two given permutations with the same parity belong to different class.

\begin{theorem}\label{ThmCl2}
Let $\mathrm{par} (\lambda)=\mathrm{par} (\mu) = (a,b)$. The following criteria characterise the conditions when $\lambda$ and $\mu$ are in the same class.
\begin{enumerate}[label=(\arabic*)]
    \item If $a=b=n$, then $\mathrm{class}(\lambda)=\mathrm{class}(\mu)$.
    \item If $a=b=n-1$, then $\mathrm{class}(\lambda)=\mathrm{class}(\mu)$ precisely when $\{\lambda_n, \lambda_{n+1}\}=\{\mu_n,\mu_{n+1}\}.$
    \item If $1<a,b\le n-1$, $(a,b)\neq (n-1,n-1)$, then
    $\mathrm{class}(\lambda)=\mathrm{class}(\mu)$ precisely when $\lambda$ and $\mu$ are friends.
    \item If exactly one of $a$ or $b$ equals $1$, then
    $\mathrm{class}(\lambda)=\mathrm{class}(\mu)$ precisely when there exists $\mu'\in\mathrm{fam}(\mu)$ such that $\lambda$ and $\mu'$ are friends.
    \item If $a=b=1$, then $\mathrm{class}(\lambda)=\mathrm{class}(\mu)$ precisely when $\mu\in\mathrm{fam}(\lambda).$
\end{enumerate}
\end{theorem}

With this in mind, we are ready to prove Theorem \ref{ThmCl1}.

\begin{proof}[Proof of Theorem \ref{ThmCl1}.]

This is done by simple counting using the results in Theorem \ref{ThmCl2}. If $\mathrm{par}(\sigma)=(a,b)$, then
\[
    \left|\mathrm{class}(\sigma)\right| =
    \begin{cases}
    2n!n! &\mbox{if $a=b=n$,}\\
    2(n-1)!(n-1)! &\mbox{if $a=b=n-1$,}\\
    \mathrm{max}\{a,2\}!\mathrm{max}\{b,2\}! &\mbox{otherwise.}
    \end{cases}
\]

For any integer $a,b \in [n]$, we define $\mathcal{P}(a,b)=\big|\{\sigma \,\colon\, \mathrm{par}(\sigma)=(a,b)\}\big|$.
By simple counting, we observe that if $a,b < n$,
\vspace{-3mm}
\begin{align*}
    \frac{\mathcal{P}(a,b)}{2n!n!} = \displaystyle{{2n-2-a-b}\choose{n-2}}+\displaystyle{{2n-2-a-b}\choose{n-1-a}}
\end{align*}
while $\mathcal{P}(n,n) = 2n!n!$, and $\mathcal{P}(a,b) = 0$ otherwise.
For any pair $(a,b)$, the number of different maps generated by $\sigma$ such that $\mathrm{par}(\sigma) = (a,b)$ is $\mathcal{P}(a,b)/\left|\mathrm{class}(\sigma)\right|$.
The result thus follows by taking the sum of all possible $a,b \in [n]$.
\end{proof}

\begin{remark}
    The number of different permutation-generated maps for $n = 1$ up to $n = 6$ is given by
\[1, 3, 154, 8369, 711226, 90349957.\]
This is a new sequence in the Online Encyclopedia of Integer Sequences (OEIS) \cite{sloane2003line} with sequence number \href{https://oeis.org/A344898}{A344898}. Moreover, one can verify that the ratio between the number of different maps and the number of all possible permutations asymptotically converges to $(\sqrt{e}-5/4)^2 \approx 0.158979$. 
\end{remark}

\subsection{A reformulation of a combinatorial identity}
In \cite[Section 3]{callan2009combinatorial}, Callan gave the following double factorial identity whose proof utilises the Hafnian of a particular matrix. For a given Dyck path $P$, the sets $\mathbf{U}_P = \{u_1<u_2<\ldots<u_n\,\colon\,P_{u_i}=\mathtt{u}\}$ and $\mathbf{D}_P =\{d_1<d_2<\ldots<d_n\,\colon\,P_{d_i}=\mathtt{d}\}$ encode all the up-steps and down-steps of $P$ respectively.

\begin{theorem}[Proposition 1 of \cite{callan2009combinatorial}]\label{ThmCallan}
    Let $P$ be a Dyck path of size $n$ and $h_{u_k}$ be the height of step $u_k \in \mathbf{U}_P$. We have that
    \vspace{-3mm}
    \begin{equation}\label{eq1}
    \sum_{P\in\mathcal{D}_n}\prod_{k=1}^n h_{u_k}(P)=(2n-1)!!.
\end{equation}
\end{theorem}

We will offer a different proof of Theorem \ref{ThmCallan} by defining a weaker version of the equivalence relation $\sim$ defined previously.

\begin{proof}
    Let $Q\in \mathcal{D}_n$ be fixed. We define an equivalence relation $\sim_Q$ on $S_{2n}$ with $\lambda \sim_Q \mu$ if and only if $\lambda(Q)=\mu(Q)$. Given two Dyck paths $P$ and $Q$, we proceed by finding an algorithm to find $\sigma$ satisfying $\sigma(Q) = P$.
        \begin{algo}
        Given Dyck paths $P$ and $Q$, find $\sigma$ such that $\sigma(Q)=P$.
            \begin{enumerate}
                \item Set $\mathbf{U}_P = \{u_1<u_2<\ldots<u_n\,\colon\,P_{u_i}=\mathtt{u}\}$ and $\mathbf{D}_P=\{d_1<d_2<\ldots<d_n\,\colon\,P_{d_i}=\mathtt{d}\}$.
                \item For $i\in [n]$, assign $\sigma_{u_i}$ freely such that all of $\sigma_{u_1},\sigma_{u_2},\dots ,\sigma_{u_n}$ belong to different tunnels of $Q$, that is, $\tau_Q(\sigma_{u_i})\neq \sigma_{u_j}$ for any $i\neq j$.
                \item Iterate for $k=1,\dots, n$:
                     \begin{enumerate}
                         \item Gather all indices $i$ which have the properties that $u_{i}<d_k$ and $\tau_Q(\sigma_{u_{i}})$ has not been filled yet. In other words, construct the set $\{i \,\colon\, u_i<d_k \text{ and } \tau_Q(\sigma_{u_i})\neq \sigma_{d_j}\text{ for some } j<k \}$. Choose any element $i^\ast$ from the set. 
	                    \item Assign $\sigma_{d_k}=\tau_Q(u_{i^*})$.
	                \end{enumerate}
            \end{enumerate}    
        \end{algo}
    Note that the number of choices of the assignment of $\sigma_{u_i}$ in Step 2 is $2^n n!$. Also, at iteration $k$ of Step 3, the number of ways to choose $i^*$ is the same as the height of $P$ at the beginning of step $d_k$, which is $h_{d_k-1}(P)$. Thus, the number of ways to choose $i^*$ for all iterations is 
        \vspace{-3mm}
        \begin{align*}
            L_P := \prod_{k=1}^n h_{d_k-1}(P) =\prod_{k=1}^n h_{u_k}(P).
        \end{align*}
    So, the total number of possible $\sigma$ is $2^n n! L_P$. As a result, when $Q$ is fixed, the size of each equivalence class $C_P := \{\sigma \,\colon\, \sigma(Q)=P\}$ induced by $\sim_Q$ is $|C_P|=2^n n!L_P$ for all $P$. Our derivation of the identity (\ref{eq1}) follows since $\sum_{P\in \mathcal{D}_n} |C_P|=|S_{2n}|=(2n)!$. 
\end{proof}

\section{Permutation-generated bijections}
We observe that this generalisation does not retain bijectivity. It is therefore of interest to classify all permutations whose maps they generate are bijections. Such a permutation has a nice geometric interpretation which we call the \emph{circularly-connected permutation} or CCP in short.

A set $X \subseteq [2n]$ is a \emph{block} of size $k$ if there exists some $x \in X$ such that $X = \{x,x+1, \ldots, x+k-1\}$ where the addition is regarded modulo $2n$. This is best illustrated if we order the elements of $[2n]$ on a circle in a clockwise manner.

\vspace{-3mm}
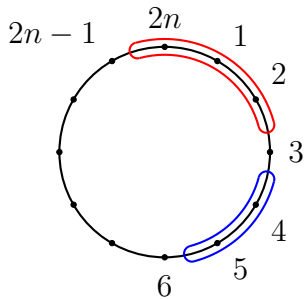
\begin{figure}[H]
    \centering
    \begin{tikzpicture}[scale=0.7]
    
        \draw[thick] (0,0) circle (2cm);
        \draw [red,thick,domain=195:375] plot ({1.93185+0.15*cos(\x)}, {0.51764+0.15*sin(\x)});
        \draw [red,thick,domain=105:285] plot ({-0.517638+0.15*cos(\x)}, {1.93185+0.15*sin(\x)});
        \draw [red,thick,domain=15:105] plot ({2.15*cos(\x)}, {2.15*sin(\x)});
        \draw [red,thick,domain=15:105] plot ({1.85*cos(\x)}, {1.85*sin(\x)});
        
        \draw [blue,thick,domain=-15:165] plot ({1.93185+0.15*cos(\x)}, {-0.51764+0.15*sin(\x)});
        \draw [blue,thick,domain=105:285] plot ({0.517638+0.15*cos(\x)}, {-1.93185+0.15*sin(\x)});
        \draw [blue,thick,domain=285:345] plot ({2.15*cos(\x)}, {2.15*sin(\x)});
        \draw [blue,thick,domain=285:345] plot ({1.85*cos(\x)}, {1.85*sin(\x)});
        \coordinate (A) at (1,1.73) {};
        \node[draw, circle, fill=black, minimum size=2pt,inner sep=0pt,radius=1pt] at (A) {};
        \node at (A) [above right = 0.866mm of A] {$1$};
        \coordinate (B) at (1.73,1) {};
        \node[draw, circle, fill=black, minimum size=2pt,inner sep=0pt,radius=1pt] at (B) {};
        \node at (B) [above right = 0.866mm of B] {$2$};
        \coordinate (C) at (2,0) {};
        \node[draw, circle, fill=black, minimum size=2pt,inner sep=0pt,radius=1pt] at (C) {};
        \node at (C) [right = 1mm of C] {$3$};
        \coordinate (D) at (1.73,-1) {};
        \node[draw, circle, fill=black, minimum size=2pt,inner sep=0pt,radius=1pt] at (D) {};
        \node at (D) [below right = 0.866mm of D] {$4$};
        \coordinate (E) at (1,-1.73) {};
        \node[draw, circle, fill=black, minimum size=2pt,inner sep=0pt,radius=1pt] at (E) {};
        \node at (E) [below right = 0.866mm of E] {$5$};
        \coordinate (F) at (0,-2) {};
        \node[draw, circle, fill=black, minimum size=2pt,inner sep=0pt,radius=1pt] at (F) {};
        \node at (F) [below = 1mm of F] {$6$};
        \coordinate (G) at (-1,-1.73) {};
        \node[draw, circle, fill=black, minimum size=2pt,inner sep=0pt,radius=1pt] at (G) {};
        \coordinate (H) at (-1.73,-1) {};
        \node[draw, circle, fill=black, minimum size=2pt,inner sep=0pt,radius=1pt] at (H) {};
        \coordinate (I) at (-2,0) {};
        \node[draw, circle, fill=black, minimum size=2pt,inner sep=0pt,radius=1pt] at (I) {};
        \coordinate (J) at (-1.73,1) {};
        \node[draw, circle, fill=black, minimum size=2pt,inner sep=0pt,radius=1pt] at (J) {};
        \coordinate (K) at (-1,1.73) {};
        \node[draw, circle, fill=black, minimum size=2pt,inner sep=0pt,radius=1pt] at (K) {};
        \node at (K) [above left = 0.866mm of K] {$2n-1$};
        \coordinate (L) at (0,2) {};
        \node[draw, circle, fill=black, minimum size=2pt,inner sep=0pt,radius=1pt] at (L) {};
        \node at (L) [above = 1mm of L] {$2n$};
        
    \end{tikzpicture}
    \caption{Two blocks of size $2$ and $3$.}
    \label{fig:component}
\end{figure}

\vspace{-5mm}
\begin{defn}
A permutation $\sigma \in S_{2n}$ is a CCP if for each $k \in [2n]$, $\sigma_{[k]}$ is a block of size $k$. The set of all CCPs on $[2n]$ is denoted by $\mathrm{CCP}_{2n}$.
\end{defn}

\begin{exmp}
The permutation $\sigma = 213645 \in S_6$ is a CCP while $\sigma' = 236145$ is not.
\end{exmp}

\begin{remark}
    Using the language of graph theory, a permutation $\sigma$ is a CCP if the subgraph of $C_{2n}$ induced by $\sigma_{[k]}$ is connected for all $k\in [2n]$.
\end{remark}

We proceed by establishing the main claim that CCPs generate bijections.

\begin{theorem}[Characterisation of permutation-generated bijections]\label{ThmCCP}
    $\sigma(\cdot)$ is a bijection if and only if $\sigma$ is a CCP.
\end{theorem}
We divide the discussion of the proof into two subsections, one for each direction.

\subsection{CCPs generate bijections}

We show that we can reverse the construction of $\sigma$-path given in Definition \ref{Def1}. In other words, given $\sigma$ and a Dyck path $Q$, we can find $P$ such that $\sigma(P) = Q$.

\begin{defn}
Let $\tau$ be a \emph{pairing permutation} of $[2n]$, that is, an involution with no fixed points. For any $a,b\in [2n]$, we say that the tuple $(a,b)$ is $\tau$-\emph{non-crossing} if for any block $C$ with endpoints $a$ and $b$, $\tau(c)\in C$ for every $c\in C$. The permutation $\tau$ is called \emph{non-crossing} if $(c,\tau(c))$ is $\tau$-non-crossing for all $c \in [2n]$.
\end{defn}

It follows that $\tau$ is the tunneling of a Dyck path if and only if $\tau$ is non-crossing. Note that $(a,b)$ being $\tau$-non-crossing is equivalent with $(a,b)$ being a tunnel pair in a Dyck path whose tunneling is $\tau$. We are now ready to give an inverse algorithm to find $P$ given $\sigma$ and $Q$ such that $\sigma(P) = Q$.

\begin{algo}[Inverse Algorithm]\label{AlgInvAlgebraic} Given a Dyck path $Q$ and $\sigma \in \mathrm{CCP}_{2n}$, find a Dyck path $P$ and $\tau \in S_{2n}$ such that $\sigma(P) = Q$ and $\tau_P = \tau$.

Consider the circular representation of $Q$ and initially set all elements of $[2n]$ to be unpaired. For all $k\in [2n]$ such that $Q_k = \mathtt{d}$, do the following procedure:
\begin{enumerate}
    \item Find an element $w$ such that $w$ is an endpoint of the block $\sigma_{[k-1]}$ and $w$ is next to $\sigma_k$.
    \item From $w$, traverse the block $\sigma_{[k-1]}$ to find the first unpaired element $v \in \sigma_{[k-1]}$.
    \item Set $\tau(\sigma_k)=v$ and $\tau(v)=\sigma_k$.
    \item Set both $\sigma_k$ and $v$ to be paired.
\end{enumerate}
Once all iterations are done, set $P$ to be the Dyck path whose tunneling is $\tau$.
\end{algo}

\begin{exmp}
Set $Q = \mathtt{uududd}$ and $\sigma = 162354 \in \mathrm{CCP}_6$. Note that $Q_3 = Q_5 = Q_6 = \mathtt{d}$. We illustrate Algorithm \ref{AlgInvAlgebraic} in the following table.
\begin{table}[H]
    \begin{center}
    \begin{tabular}{|c|c|c|c|c|c|c|}
        \hline
        Iteration & $k$ & $\sigma_{[k-1]}$ & $\sigma_k$ & $w$ & $v$ & Tunneling produced\\
        \hline
        1 & 3 & $\{1,6\}$ & 2 & 1 & 1 & $\tau(2)=1, \tau(1)=2$\\
        2 & 5 & $\{1,6,2,3\}$ & 5 & 6 & 6 & $\tau(5)=6, \tau(6)=5$ \\
        3 & 6 & $\{1,6,2,3,5\}$ & 4 & 3 or 5 & 3 & $\tau(3)=4, \tau(4)=3$ \\
        \hline
    \end{tabular}
    \end{center}
    \caption{Iteration table of Algorithm \ref{AlgInvAlgebraic}.}
\end{table}
The Dyck path $P$ whose tunneling is $\tau = 214365$ is $\mathtt{ududud}$.
\end{exmp}

\begin{figure}[H]
    \centering
    \begin{tikzpicture}[scale=1]
    
        \draw[thick] (0,0) circle (1cm);
        \draw[thick] (-6,0) circle (1cm);
        \draw[thick] (6,0) circle (1cm);
        
        \coordinate (A) at (-5.5,0.865);
        \node[draw, circle, fill=black, minimum size=2pt,inner sep=0pt,radius=1pt] at (A) {};
        \node at (A) [above right = 0.866mm of A] {$1 = w = v$};
        \coordinate (B) at (-5,0);
        \node[draw, circle, fill=black, minimum size=2pt,inner sep=0pt,radius=1pt] at (B) {};
        \node at (B) [right = 1mm of B] {$2 = \sigma_3$};
        \coordinate (C) at (-5.5,-0.865);
        \node[draw, circle, fill=black, minimum size=2pt,inner sep=0pt,radius=1pt] at (C) {};
        \node at (C) [below right = 0.866mm of C] {$3$};
        \coordinate (D) at (-6.5,-0.865);
        \node[draw, circle, fill=black, minimum size=2pt,inner sep=0pt,radius=1pt] at (D) {};
        \node at (D) [below left = 0.866mm of D] {$4$};
        \coordinate (E) at (-7,0);
        \node[draw, circle, fill=black, minimum size=2pt,inner sep=0pt,radius=1pt] at (E) {};
        \node at (E) [left = 1mm of E] {$5$};
        \coordinate (F) at (-6.5,0.865);
        \node[draw, circle, fill=black, minimum size=2pt,inner sep=0pt,radius=1pt] at (F) {};
        \node at (F) [above left = 0.866mm of F] {$6$};
        
        \coordinate (G) at (0.5,0.865);
        \node[draw, circle, fill=black, minimum size=2pt,inner sep=0pt,radius=1pt] at (G) {};
        \node at (G) [above right = 0.866mm of G] {$1$};
        \coordinate (H) at (1,0);
        \node[draw, circle, fill=black, minimum size=2pt,inner sep=0pt,radius=1pt] at (H) {};
        \node at (H) [right = 1mm of H] {$2$};
        \coordinate (I) at (0.5,-0.865);
        \node[draw, circle, fill=black, minimum size=2pt,inner sep=0pt,radius=1pt] at (I) {};
        \node at (I) [below right = 0.866mm of I] {$3$};
        \coordinate (J) at (-0.5,-0.865);
        \node[draw, circle, fill=black, minimum size=2pt,inner sep=0pt,radius=1pt] at (J) {};
        \node at (J) [below left = 0.866mm of J] {$4$};
        \coordinate (K) at (-1,0);
        \node[draw, circle, fill=black, minimum size=2pt,inner sep=0pt,radius=1pt] at (K) {};
        \node at (K) [left = 1mm of K] {$\sigma_5 = 5$};
        \coordinate (L) at (-0.5,0.865);
        \node[draw, circle, fill=black, minimum size=2pt,inner sep=0pt,radius=1pt] at (L) {};
        \node at (L) [above left = 0.866mm of L] {$v = w = 6$};
        
        \coordinate (M) at (6.5,0.865);
        \node[draw, circle, fill=black, minimum size=2pt,inner sep=0pt,radius=1pt] at (M) {};
        \node at (M) [above right = 0.866mm of M] {$1$};
        \coordinate (N) at (7,0);
        \node[draw, circle, fill=black, minimum size=2pt,inner sep=0pt,radius=1pt] at (N) {};
        \node at (N) [right = 1mm of N] {$2$};
        \coordinate (O) at (6.5,-0.865);
        \node[draw, circle, fill=black, minimum size=2pt,inner sep=0pt,radius=1pt] at (O) {};
        \node at (O) [below right = 0.866mm of O] {$3 = v$};
        \coordinate (P) at (5.5,-0.865);
        \node[draw, circle, fill=black, minimum size=2pt,inner sep=0pt,radius=1pt] at (P) {};
        \node at (P) [below left = 0.866mm of P] {$\sigma_6 = 4$};
        \coordinate (Q) at (5,0);
        \node[draw, circle, fill=black, minimum size=2pt,inner sep=0pt,radius=1pt] at (Q) {};
        \node at (Q) [left = 1mm of Q] {$w = 5$};
        \coordinate (R) at (5.5,0.865);
        \node[draw, circle, fill=black, minimum size=2pt,inner sep=0pt,radius=1pt] at (R) {};
        \node at (R) [above left = 0.866mm of R] {$6$};
        
        \draw [red,thick,domain=240:420] plot ({-5.5+0.15*cos(\x)}, {0.866+0.15*sin(\x)});
        \draw [red,thick,domain=120:300] plot ({-6.5+0.15*cos(\x)}, {0.866+0.15*sin(\x)});
        \draw [red,thick,domain=60:120] plot ({-6+1.15*cos(\x)}, {1.15*sin(\x)});
        \draw [red,thick,domain=60:120] plot ({-6+0.85*cos(\x)}, {0.85*sin(\x)});
        
        \draw[->,>=stealth',thick, blue] (-5.3,0) arc[radius=0.5, start angle=0, end angle=90];
        
        \draw [red,thick,domain=120:300] plot ({-0.5+0.15*cos(\x)}, {0.866+0.15*sin(\x)});
        \draw [red,thick,domain=120:300] plot ({0.5+0.15*cos(\x)}, {-0.866+0.15*sin(\x)});
        \draw [red,thick,domain=-60:120] plot ({1.15*cos(\x)}, {1.15*sin(\x)});
        \draw [red,thick,domain=-60:120] plot ({0.85*cos(\x)}, {0.85*sin(\x)});
        
        \draw[->,>=stealth',thick, blue] (-0.75,0) arc[radius=0.5, start angle=180, end angle=90];
        
        \draw [red,thick,domain=180:360] plot ({5+0.15*cos(\x)}, {0.15*sin(\x)});
        \draw [red,thick,domain=120:300] plot ({6.5+0.15*cos(\x)}, {-0.866+0.15*sin(\x)});
        \draw [red,thick,domain=-60:180] plot ({6+1.15*cos(\x)}, {1.15*sin(\x)});
        \draw [red,thick,domain=-60:180] plot ({6+0.85*cos(\x)}, {0.85*sin(\x)});
        
        \draw[->,>=stealth',thick, blue] (5.7,-0.52) arc[radius=0.6, start angle=240, end angle=-60];

    \end{tikzpicture}
    \caption{Iteration 1 (left), iteration 2 (middle), iteration 3 (right).}
\end{figure}
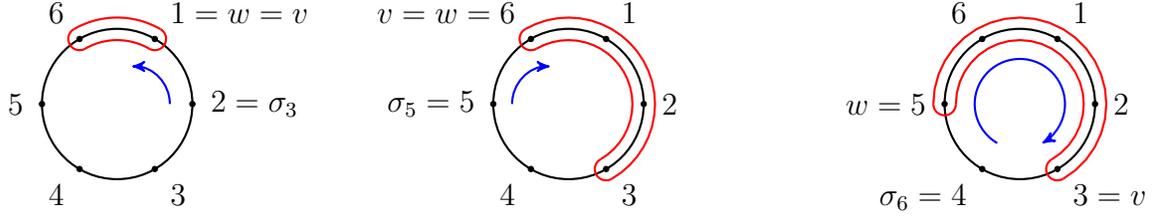

\vspace{-3mm}
Algorithm \ref{AlgInvAlgebraic} is well-defined whenever $\sigma\in \mathrm{CCP}_{2n}$. To see this, it is sufficient to guarantee the existence and uniqueness of the elements $w$ and $v$ specified in the first and second steps of each iteration (except in the last iteration) in the algorithm.
    
For any $d_i\in\mathbf{D}_Q$ with $i < n$, observe that at iteration $d_i > 1$, it is necessary that $\sigma_{d_i}$ is a consecutive element of an endpoint $w$ of the block $\sigma_{[d_i-1]}$ by the definition of CCP. Since $\sigma$ is a CCP and $i<n$, $w$ is uniquely determined.

Next, at the beginning of the second step, there are exactly $i-1$ paired elements. Since $Q$ is a Dyck path, there are at least $2i-1$ elements in $\sigma_{[k_i-1]}$. Thus, there exists at least one unpaired element in $\sigma_{[k_i-1]}$ implying the existence of $v$. Since $w$ is unique, so is $v$.

The following two lemmas establish that the output $P$ of Algorithm \ref{AlgInvAlgebraic} does indeed satisfy $\sigma(P) = Q$.

\begin{lemma}
    The permutation $\tau$ produced in Algorithm \ref{AlgInvAlgebraic} is non-crossing. Hence, $\tau$ is a tunneling of a Dyck path $P$.
\end{lemma}
\begin{proof}
    It is obvious that $\tau$ is a pairing permutation. For any $d_i\in\mathbf{D}_Q$, we claim that $(\sigma_{d_i},\tau(\sigma_{d_i}))$ is $\tau$-non-crossing. Consider the block $C$ with endpoints $\sigma_{d_i}$ and $\tau(\sigma_{d_i})$ which is also a subset of $\sigma_{[d_i]}$. Note that all elements in $C$ are paired at the end of iteration $i$.
    
    For contradiction, suppose that there exists $c\in C$ such that $\tau(c)\notin C$. Let $c=\sigma_{d_\ell}$ with $\ell<i$. At the end of iteration $d_\ell$, by the definition of the algorithm, all the elements in the block with endpoints $c$ and $\tau(c)$ which also a subset of $\sigma_{[d_\ell]}$ are paired. Observe that $\tau(\sigma_{d_i})$ must be in the block since otherwise $\tau(c)$ would be in $C$. However, $\tau(\sigma_{d_i})$ is unpaired during this iteration, a contradiction. The result follows from the claim.
\end{proof}

\begin{lemma}[Algorithm \ref{AlgInvAlgebraic} is an inverse algorithm]\label{Lem2}
    Let $\sigma \in \mathrm{CCP}_{2n}$ and $Q\in \mathcal{D}_n$. If $P$ is the Dyck path produced by Algorithm \ref{AlgInvAlgebraic}, then $\sigma(P)=Q$. 
\end{lemma}
\begin{proof}
    It is sufficient to prove that $\sigma(P)_{u_i}=Q_{u_i}=\mathtt{u}$ for all $u_i\in \mathbf{U}_Q$. From the definition of the algorithm, we have $\tau_P(\sigma_{u_j})\neq\sigma_{u_i}$ for all $u_j<u_i$. For any $k$ such that $d_k<u_i$, it is necessary that $\tau_P(\sigma_{d_k})\neq \sigma_{u_i}$ since $\sigma_{d_k}$ is paired at iteration $k$ and $\sigma_{u_i} \notin \sigma_{[d_k]}$. Therefore, $\tau_P(\sigma_\ell)\neq \sigma_{u_i}$ for all $\ell<u_i$, which means $\sigma(P)_{u_i}=Q_{u_i}=\mathtt{u}$, hence the result. 
\end{proof}

\begin{proof}[Proof of Theorem \ref{ThmCCP}($\Longleftarrow$)]
    By Lemma \ref{Lem2}, $\sigma(\cdot)$ is a surjection between finite sets.
\end{proof}
\subsection{Non-CCPs generate non-bijections}
We show that if $\sigma$ is not a CCP, there exists a Dyck path $Q$ such that $Q\not\in \mathrm{Image\,}\sigma(\cdot)$.

\begin{proof}[Proof of Theorem \ref{ThmCCP}($\Longrightarrow$)]
By contradiction, suppose that $\sigma(\cdot)$ is a bijection and $\sigma$ is not a CCP. Since $\sigma$ is not a CCP, there exists the smallest integer $k\in[2, 2n-2]$ such that $\sigma_{[k-1]}$ is a block but $\sigma_{[k]}$ is not. We consider two cases: 
    \paragraph{Case 1: $k = 2j$} 
    Consider two different blocks $B_1, B_2$, both with endpoints $\sigma_k$ and $\sigma_{k-1}$. Since $\sigma_{[k-1]}$ is a block but $\sigma_{[k]}$ is not, we have $\sigma_{[k]} \cap B_1 \neq B_1$ and $\sigma_{[k]} \cap B_2 \neq B_2$. Without loss of generality, suppose that $\sigma_{k+1}\in B_1$. Since $\sigma_{[k]} \cap B_2 \neq B_2$, there exists the smallest integer $l>k+1$ such that $\sigma_l\in B_2$.
    
    We can choose $Q$ such that $Q_{1}Q_2\ldots Q_k = (\texttt{ud})^j$ and 
    \[
    Q_{k+1}\dots Q_l =\begin{cases}
        \texttt{u}^{l-k-1}\texttt{d} & \text{ if } l\leq n+j+1,\\
        \texttt{u}^{n-j}\texttt{d}^{l-n-j} & \text{ if } l>n+j+1.
    \end{cases}
    \]
    Since $\sigma(\cdot)$ is a bijection, there exists $P\in \mathcal{D}_n$ such that $\sigma(P)=Q$. Since $\sigma(P)_l=Q_l=\mathtt{d}$ and $Q_1\ldots Q_k\in \mathcal{D}_j$, we have $\tau_P(\sigma_l)\in\{\sigma_{k+1},\dots,\sigma_{l-1}\}$. Thus, $\sigma_l\in B_2$ but $\tau_P(\sigma_l)\in B_1$. In other words, $(\sigma_k,\sigma_{k-1})=(\sigma_k,\tau_P(\sigma_k))$ is not $\tau_P$-non-crossing, a contradiction.
    
    \paragraph{Case 2: $k = 2j+1$}
    First, notice that there exists a Dyck path $D$ such that $\sigma(D)=(\mathtt{ud})^n$ since $\sigma(\cdot)$ is a bijection. Thus, $\tau_D(\sigma_{2i-1})=\sigma_{2i}$ for all $i \in [n]$, so $\sigma_{2i-1} \not \equiv \sigma_{2i} \,\mathrm{mod}\, 2$ for all $i \in [n]$.
    
    We claim that for any $Q,P\in\mathcal{D}_n$ with $Q_1\ldots Q_k=\mathtt{u}(\mathtt{ud})^j$ and $\sigma(P)=Q$, the values of $\tau_P(\sigma_3),\ldots,\tau_P(\sigma_{2j+1})$ do not depend on the choice of $Q_{k+1}\dots Q_{2n}$. We proceed with induction. For $\sigma_3$, since $\sigma(P)_3=Q_3=\mathtt{d}$, we have $\tau_P(\sigma_3)\in\{\sigma_1,\sigma_2\}$. Since $\sigma_1\not\equiv \sigma_2 \,\mathrm{mod}\,2$, there is exactly one possible option for $\tau_P(\sigma_3)$, which depends only on $\sigma$. Now, suppose that $\tau_P(\sigma_3),\tau_P(\sigma_5),\ldots,\tau_P(\sigma_{2q-1})$ satisfy the claim for some $q\geq 2$. Since $\sigma(P)_{2i-1}=Q_{2i-1}=\mathtt{d}$ for $i\in [3, q]$, we have $\tau_P(\sigma_{2i-1})\in\sigma_{[2q-1]}$. Notice that there exists exactly one element $\sigma_r\in \{\sigma_1\}\cup \{\sigma_{2i}\}_{i=1}^{q-1}$ other than $\tau_P(\sigma_3),\tau_P(\sigma_5),\ldots,\tau_P(\sigma_{2q-1})$. By the induction hypothesis, the value of $r$ does not depend on $Q_{k+1}\ldots Q_{2n}$. Because $\sigma_{2i-1} \not \equiv \sigma_{2i} \,\mathrm{mod}\, 2$ for all $ i\in [n]$, $\sigma_{[2q]}$ consists of exactly $q$ odd and $q$ even integers. Excluding all $\sigma_3,\dots,\sigma_{2q-1}$ and $\tau_P(\sigma_3),\tau_P(\sigma_5),\ldots,\tau_P(\sigma_{2q-1})$, we have  $\sigma_r \not \equiv \sigma_{2q} \,\mathrm{mod}\, 2$. Since $\sigma(P)_{2q+1}=Q_{2q+1}=\mathtt{d}$, we have $\tau_P(\sigma_{2q+1})\in\{\sigma_r,\sigma_{2q}\}$. Since $\sigma_r\not\equiv \sigma_{2q} \,\mathrm{mod}\,2$, there is exactly one possible option for $\tau_P(\sigma_{2q+1})$ that depends only on $\sigma$, which completes the induction and proving the claim.
        
    From the claim, there exists an integer $m<2j+1$ such that
    $\tau_P(\sigma_{2j+1})=\sigma_m$ for all $Q, P\in\mathcal{D}_n$ with $Q_1\dots Q_k=\mathtt{u}(\mathtt{ud})^j$ and $\sigma(P)=Q$. Consider an integer $r$ where $\sigma_r\in \{\sigma_1\}\cup \{\sigma_{2i}\}_{i=1}^{j}$ other than $\tau_P(\sigma_3),\tau_P(\sigma_5),\ldots,\tau_P(\sigma_k)$. Consider two different blocks $B_1, B_2$, both with endpoints $\sigma_k$ and $\sigma_m$. Since $\sigma_{[k-1]}$ is a block but $\sigma_{[k]}$ is not a block, we have $\sigma_{[k]} \cap B_1 \neq B_1$ and $\sigma_{[k]} \cap B_2 \neq B_2$. Without loss of generality, suppose that $\sigma_r\in B_1$. Since $\sigma_{[k]} \cap B_2 \neq B_2$, there exists the smallest integer $l\geq k+1$ such that $\sigma_l\in B_2$.
    
    Now, consider the Dyck path $Q'$ such that $Q'_{1}Q'_2\ldots Q'_k =\texttt{u} (\texttt{ud})^j$ and 
    \[
    Q'_{k+1}\dots Q'_l =\begin{cases}
        \texttt{d} & \text{ if } l=k+1,\\ 
        \texttt{u}^{l-k-1}\texttt{d} & \text{ if } k+1<l\leq n+j,\\
        \texttt{u}^{n-j-1}\texttt{d}^{l-n-j+1} & \text{ if } l>n+j.
    \end{cases}
    \]
    Since $\sigma(\cdot)$ is a bijection, there exists $P'\in \mathcal{D}_n$ such that $\sigma(P')=Q'$. Note that values of $m, r, l$ defined before also hold for $Q'$ and $P'$. Since $\sigma(P')_l=Q'_l=\mathtt{d}$ and $Q'_1\ldots Q'_k=\mathtt{u}(\mathtt{ud})^j$, we have $\tau_{P'}(\sigma_l)\in\{\sigma_r\}\cup\{\sigma_i\}_{i=k+1}^{l-1}$. Thus, $\sigma_l\in B_2$ but $\tau_{P'}(\sigma_l)\in B_1$. In other words, $(\sigma_k,\sigma_m)=(\sigma_k,\tau_{P'}(\sigma_k))$ is not $\tau_{P'}$-non-crossing, a contradiction. See the figure below for illustration.
    
    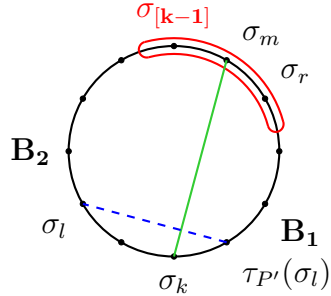
\begin{figure}[H]\label{Fig5}
    \centering
    \begin{tikzpicture}[scale=0.7]
    
        \draw[thick] (0,0) circle (2cm);
        \draw [red,thick,domain=195:375] plot ({1.93185+0.15*cos(\x)}, {0.51764+0.15*sin(\x)});
        \draw [red,thick,domain=105:285] plot ({-0.517638+0.15*cos(\x)}, {1.93185+0.15*sin(\x)});
        \draw [red,thick,domain=15:105] plot ({2.15*cos(\x)}, {2.15*sin(\x)});
        \draw [red,thick,domain=15:105] plot ({1.85*cos(\x)}, {1.85*sin(\x)});
        
        \coordinate (A) at (1,1.73) {};
        \node[draw, circle, fill=black, minimum size=2pt,inner sep=0pt,radius=1pt] at (A) {};
        \node at (A) [above right = 0.866mm of A] {$\sigma_m$};
        \coordinate (B) at (1.73,1) {};
        \node[draw, circle, fill=black, minimum size=2pt,inner sep=0pt,radius=1pt] at (B) {};
        \node at (B) [above right = 0.866mm of B] {$\sigma_r$};
        \coordinate (C) at (2,0) {};
        \node[draw, circle, fill=black, minimum size=2pt,inner sep=0pt,radius=1pt] at (C) {};
        \node at (C) [right = 1mm of C] {};
        \coordinate (D) at (1.73,-1) {};
        \node[draw, circle, fill=black, minimum size=2pt,inner sep=0pt,radius=1pt] at (D) {};
        \node at (D) [below right = 0.866mm of D] {$\mathbf{B_1}$};
        \coordinate (E) at (1,-1.73) {};
        \node[draw, circle, fill=black, minimum size=2pt,inner sep=0pt,radius=1pt] at (E) {};
        \node at (E) [below right = 0.866mm of E] {$\tau_{P'}(\sigma_l)$};
        \coordinate (F) at (0,-2) {};
        \node[draw, circle, fill=black, minimum size=2pt,inner sep=0pt,radius=1pt] at (F) {};
        \node at (F) [below = 1mm of F] {$\sigma_k$};
        \coordinate (G) at (-1,-1.73) {};
        \node[draw, circle, fill=black, minimum size=2pt,inner sep=0pt,radius=1pt] at (G) {};
        \coordinate (H) at (-1.73,-1) {};
        \node[draw, circle, fill=black, minimum size=2pt,inner sep=0pt,radius=1pt] at (H) {};
        \node at (H) [below left = 0.866mm of H] {$\sigma_l$};
        \coordinate (I) at (-2,0) {};
        \node[draw, circle, fill=black, minimum size=2pt,inner sep=0pt,radius=1pt] at (I) {};
        \node at (I) [left = 1mm of I] {$\mathbf{B_2}$};
        \coordinate (J) at (-1.73,1) {};
        \node[draw, circle, fill=black, minimum size=2pt,inner sep=0pt,radius=1pt] at (J) {};
        \coordinate (K) at (-1,1.73) {};
        \node[draw, circle, fill=black, minimum size=2pt,inner sep=0pt,radius=1pt] at (K) {};
        \node at (K) [above left = 0.866mm of K] {};
        \coordinate (L) at (0,2) {};
        \node[draw, circle, fill=black, minimum size=2pt,inner sep=0pt,radius=1pt] at (L) {};
        \node at (L) [above = 1mm of L] {\color{red}$\mathbf{\sigma_{[k-1]}}$};
        \draw [limegreen, thick] (A)--(F);
        \draw [blue, thick, dashed] (E)--(H);
    \end{tikzpicture}
    \caption{$(\sigma_k,\sigma_m)$ is not $\tau_{P'}$-non-crossing.}
    \end{figure}

    \vspace{-3mm}
    Therefore, $\sigma$ must be a CCP and this concludes the proof.
\end{proof}

\begin{remark} 
    Deutsch and Elizalde's permutations $\sigma$ and $\sigma^{(r)}$ in \cite{elizalde2003simple} are CCPs. In total, there are exactly $n2^{2n-1}$ CCPs. If $\sigma \in \mathrm{CCP}_{2n}, n\geq 3$, then $\mathrm{par}(\sigma)=(1,1)$. Thus $\left|\mathrm{class}(\sigma)\right|=4$ and there are exactly $n2^{2n-3}$ different permutation-generated bijections.
\end{remark}

\section{Statistics preserved by the bijections}

In this section, we define several statistics which are translated by the permutation-generated bijections to some classic statistics in $\mathcal{D}_n$.

\begin{defn}
Let $P\in \mathcal{D}_n$ and $a,k\in [2n]$. Suppose that $S\subseteq [2n]$. A step $i\in S$ is \textit{unmatched} in $S$ if $\tau_P(i)\notin S$ and \textit{matched} otherwise. Let $u_{a,k}(P)$ denote the number of unmatched steps in $\{P_i\}_{i\in I}$ where $I=\{a,\ldots, a+k-1\}\subseteq [2n]$, that is, $k$ circularly consecutive steps of $P$ starting from step $a$. 
\end{defn}

\begin{exmp}
In the figure below, unmatched and matched steps of $P$ with respect to our choice of $a$ and $k$ are represented by dashed red and thick blue lines respectively. 
\begin{figure}[H]
    \begin{center}
    \begin{tikzpicture}[scale=0.58]
    \tikzstyle{every node}=[draw,circle,fill=black,minimum size=3.5pt,inner sep=0pt,radius=7pt]
      \node (A) at (0,0) {};
      \node (B) at (1,1) {};
      \node (C) at (2,2) {};
      \node (D) at (3,3) {};
      \node (E) at (4,2) {};
      \node (F) at (5,3) {};
      \node (G) at (6,4) {};
      \node (H) at (7,3) {};
      \node (I) at (8,2) {};
      \node (J) at (9,1) {};
      \node (K) at (10,2) {};
      \node (L) at (11,1) {};
      \node (M) at (12,0) {};
      \draw (A)--(B);
      \draw [red,dashed,thick](B)--(C);
      \draw [blue,line width=0.7mm](C)--(D);
      \draw [blue,line width=0.7mm](D)--(E);
      \draw [red,dashed,thick](E)--(F);
      \draw [red,dashed,thick](F)--(G);
      \draw (G)--(H);
      \draw (H)--(I);
      \draw (I)--(J);
      \draw (J)--(K);
      \draw (K)--(L);
      \draw (L)--(M);
      
      \node (AA) at (16,0) {};
      \node (AB) at (17,1) {};
      \node (AC) at (18,2) {};
      \node (AD) at (19,3) {};
      \node (AE) at (20,2) {};
      \node (AF) at (21,3) {};
      \node (AG) at (22,4) {};
      \node (AH) at (23,3) {};
      \node (AI) at (24,2) {};
      \node (AJ) at (25,1) {};
      \node (AK) at (26,2) {};
      \node (AL) at (27,1) {};
      \node (AM) at (28,0) {};
      \draw [blue,line width=0.7mm](AA)--(AB);
      \draw [blue,line width=0.7mm](AB)--(AC);
      \draw [red, thick, dashed](AC)--(AD);
      \draw (AD)--(AE);
      \draw (AE)--(AF);
      \draw (AF)--(AG);
      \draw (AG)--(AH);
      \draw (AH)--(AI);
      \draw [blue,line width=0.7mm](AI)--(AJ);
      \draw [blue,line width=0.7mm](AJ)--(AK);
      \draw [blue,line width=0.7mm](AK)--(AL);
      \draw [blue,line width=0.7mm](AL)--(AM);
    \end{tikzpicture}
    \end{center}
    \caption{$P=\mathtt{uuuduudddudd}$, $u_{2,5}(P)=3$ (left) and $u_{9,7}(P)=1$ (right). }
\end{figure}
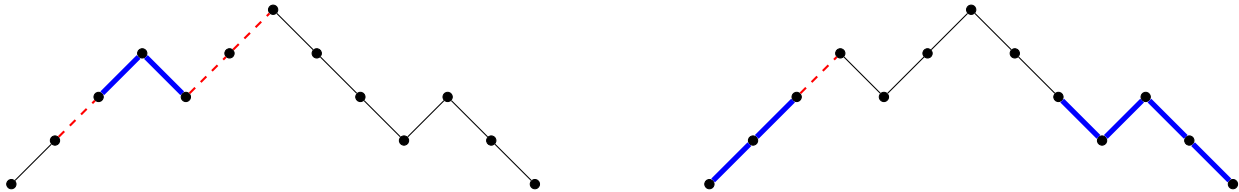

\end{exmp}

\vspace{-3mm}
It is easy to see that for any $P\in \mathcal{D}_n$ and $k\in [2n]$, we obtain $h_k(P) = u_{1,k}(P)$. For other choices of $a$, the relationship between $h_k(P)$ and $u_{a,k}(P)$ is less obvious. However, there is a nice relationship between them when $P$ is mapped by $\sigma(\cdot)$ for some $\sigma$.

\begin{theorem}\label{StatsLemma}
    Let $a, k\in [2n]$ be arbitrary. There exists a permutation $\sigma \in S_{2n}$ that depends on $a$ and $k$ such that for any $P\in \mathcal{D}_n$, we have that $u_{a,k}(P) = h_k\left(\sigma(P)\right)$.
\end{theorem}
\begin{proof}
    We specify $\sigma$ to be any CCP such that $\sigma_{[k]}=\{a,\ldots, a+k-1\}$. By the definition of $\sigma(\cdot)$, the number of down-steps among the first $k$ steps of $\sigma(P)$ is equal to the number of different tunnel pairs whose steps are taken from the set $\{P_a,\ldots ,P_{a+k-1}\}$. To simplify matters, if we denote the latter by $t$, consequently $h_k\left(\sigma(P)\right)=k-2t$. By definition, this equals to $u_{a,k}(P)$.
\end{proof}

The theorem above implies that the statistics $u_{a,k}$ and $h_k$ are distributed identically.

\begin{corollary}
    Let $a,k\in [2n]$ and $\ell \in [n]$. We have the following enumeration identity:
    \begin{equation}\label{eq2}
    \Big|\{ P\in \mathcal{D}_n: u_{a,k}(P)=\ell\}\Big| = \Big|\{ P\in \mathcal{D}_n: h_{k}(P)=\ell\}\Big|.
    \end{equation}
\end{corollary}
    
\begin{proof}
    The result follows from the bijectivity of $\sigma(\cdot)$ defined in Theorem \ref{StatsLemma}.
\end{proof}
Note that \eqref{eq2} does not depend on $a$. This means that the number of Dyck paths with certain number of unmatched steps among any $k$ consecutive (circular) steps is always the same.

\begin{remark}
    The right-hand side of \eqref{eq2} has a closed form
    \[
    \vspace{-3mm}
    \dfrac{(\ell +1)^2}{(k+1)(2n-k+1)}{k+1 \choose \frac{k+\ell}{2}+1}{2n-k+1 \choose n-\frac{k+\ell}{2}}
    \]
    with the convention that $\binom{n}{i} = 0$ if $i\notin\mathbb{Z}$.
    The value above can be easily derived by considering the condition $h_k(P)=\ell$ as a path $P$ that is forced to pass through a particular point.
\end{remark}

Theorem \ref{StatsLemma} also implies another correspondence between two other statistics in $\mathcal{D}_n$. For any $P\in \mathcal{D}_n$ and $a\in [2n]$, let $h(P)$ be the height of the highest peak in $P$.  The distribution of the number of Dyck paths with statistic $h$ is well-known, for example one could consult OEIS \href{https://oeis.org/A080936}{A080936} \cite{sloane2003line}. Obviously, $h(P)=\max\limits_{k\in[2n]}h_k(P)$. Similarly, define \[
\vspace{-3mm}
u_{\max}^{(a)}(P):= \max\limits_{k\in [2n]} u_{a,k}(P).\]
The following corollary shows that $u_{\max}^{(a)}$ and $h$ are equidistributed.
\begin{corollary}\label{lastcor}
    Let $P\in \mathcal{D}_n$ and $a\in [2n]$ be arbitrary, then for all $\ell\in n$,
    \begin{equation} \label{eq3}
    \Big|\{ P\in \mathcal{D}_n: u_{\max}^{(a)}(P)=\ell\}\Big| = \Big|\{ P\in \mathcal{D}_n: h(P)=\ell\}\Big|.
    \end{equation}
\end{corollary}
\begin{proof}
    Consider the permutation $\sigma \in S_{2n}$ with $\sigma_k=a+k-1$. Thus $\sigma$ satisfies the condition mentioned in the proof of Theorem \ref{StatsLemma} for any $k\in [2n]$. The result follows directly from this.
\end{proof}

\begin{exmp}
To illustrate Corollary \ref{lastcor}, we list the the values of $u_{\max}^{(a)}$ for $n=3$ below.
\begin{table}[H]
    \begin{center}
    \begin{tabular}{|c|c|c|c|c|c|c|c|}
        \hline
        $P$ & $h(P)$ & $u_{\max}^{(1)}(P)$ & $u_{\max}^{(2)}(P)$ & $u_{\max}^{(3)}(P)$ & $u_{\max}^{(4)}(P)$ & $u_{\max}^{(5)}(P)$ & $u_{\max}^{(6)}(P)$\\
        \hline
        $\mathtt{uuuddd}$ & 3 & 3 & 2 & 2 & 3 & 2 & 2\\
        $\mathtt{uududd}$ & 2 & 2 & 1 & 2 & 1 & 2 & 1\\
        $\mathtt{uuddud}$ & 2 & 2 & 2 & 3 & 2 & 2 & 3\\
        $\mathtt{uduudd}$ & 2 & 2 & 3 & 2 & 2 & 3 & 2\\
        $\mathtt{ududud}$ & 1 & 1 & 2 & 1 & 2 & 1 & 2\\
        \hline
    \end{tabular}
    \end{center}
    \caption{Various values of $u_{\max}^{(a)}$ in $\mathcal{D}_3$.}
\end{table}

\end{exmp}

\section*{Concluding Remarks}
Although the generalisation produces considerably large amount of bijections compared to those of Deutsch and Elizalde, the application is known for only relatively small numbers of CCPs. There might also be some relation of another subclass of CCPs to some established Dyck paths statistics. Finally, one could investigate the application of the permutation-generated maps in general by utilising the characterisation of permutations admitting the same map provided in Theorem \ref{ThmCl2}.

\section*{Acknowledgements}
We would like to thank Yohanes Tjandrawidjaja for his insight in the computational
aspect of this paper. We also thank the referees for their valuable comments.

\section*{Appendix A: Proof of Theorem \ref{ThmCl2}}

\begin{proof}[Proof of Theorem \ref{ThmCl2} (1)]
    For any $\lambda,\mu \in S_{2n}$ with $\mathrm{par} (\lambda)=\mathrm{par} (\mu) = (n,n)$ and $D \in \mathcal{D}_n$, we have $\lambda(D)_i=\mu(D)_i=\mathtt{u}$ for all $i \in [n]$ and $\lambda(D)_j=\mu(D)_j=\mathtt{d}$ for all $j \in [n+1,2n]$, hence the result.
\end{proof}

For our proof of Theorem \ref{ThmCl2} parts $(2)-(5)$, we need to prove the following lemma. In the following proof, we shall call any triple $(i,j,k)$ that satisfies Proposition \ref{Prop3} as a \emph{destroying} triple.

\begin{lemma}\label{LemmaCl3}
Let $\lambda,\mu \in S_{2n}$ with $\mathrm{par} (\lambda)=\mathrm{par} (\mu)=(a,b)$, where $a,b<n$ and $(a,b)\neq (n-1,n-1)$. If $\mathrm{class}(\lambda)=\mathrm{class}(\mu)$, then the following statements hold:
\begin{enumerate}[label=(\arabic*)]
    \item $\lambda_i=\mu_i$ for all $i \in [a+2, 2n-1-b]$,
    \item If $1<a,b<n-1$, then $\lambda_{a+1}=\mu_{a+1}$ and $\lambda_{2n-b}=\mu_{2n-b}$,
    \item If $a=n-1$ (resp. $b=n-1$), then $\lambda_{a+1}=\mu_{a+1}$ (resp. $\lambda_{2n-b}=\mu_{2n-b}$),
    \item If $a=1$ (resp. $b=1$), then $\{\lambda_1,\lambda_2\}=\{\mu_1,\mu_2\}$ $\Big($resp.  $\{\lambda_{2n-1},\lambda_{2n}\}=\{\mu_{2n-1},\mu_{2n}\} \Big)$.
\end{enumerate}
\end{lemma}

\begin{proof}
    We will prove the results by proving the contrapositive of each statement. For each statement, if there exists such a destroying triple, then we are done by Proposition \ref{Prop3}. Suppose otherwise.
    \begin{enumerate}[label=(\arabic*)]
        \item Let $\lambda_j\neq \mu_j$ for some $j \in [a+2, 2n-1-b]$. We break into two cases.
        \paragraph{Case 1: $\lambda_j \equiv \mu_j \,\mathrm{mod}\, 2$} ~\\ Let $i$ be an integer such that $\lambda_i=\lambda_j-1$. Without loss of generality, suppose that $i<j$. Let $k$ be any integer such that $j<k\leq 2n$ and $\mu_j \not \equiv \mu_k \,\mathrm{mod}\, 2$. Since $(i,j,k)$ is not a destroying triple, we must have $\mu_k=\lambda_j-1$. Let $p$ be any integer such that $\mu_j \not \equiv \mu_p \,\mathrm{mod}\, 2$ with $\mu_p\neq\lambda_j-1$ and $\mu_p\neq \lambda_j+1$. Notice that $p$ always exists for $n\geq 3$ and $p<j$ because $\mu_p \neq \lambda_j-1$. Let $q$ be an integer that satisfies $\lambda_q=\lambda_j+1$. If $j<q$, then we have a destroying triple $(p,j,q)$, a contradiction. If $j>q$, then we have a destroying triple $(q,j,k)$, a contradiction.
        
        \paragraph{Case 2: $\lambda_j \not \equiv \mu_j \,\mathrm{mod}\, 2$} Let $s$ be any integer such that $\lambda_s=\lambda_j-1$ or $\lambda_s=\lambda_j+1$ which also satisfies $\lambda_s\neq \mu_j$. Let $t$ be any integer such that $\mu_t=\mu_j-1$ or $\mu_t=\mu_j+1$ satisfying $\mu_t\neq \lambda_j$. Without loss of generality, suppose that $s<j$. If $t>j$, then $(s,j,t)$ is a destroying triple, a contradiction. If $t<j$, we consider three subcases: 
        \begin{itemize}
            \item If $\lambda_j \not \equiv \lambda_{2n-b} \,\mathrm{mod}\, 2$ and $\mu_j \not \equiv \mu_{2n-b} \,\mathrm{mod}\, 2$, consider the triple $(s,j,2n-b)$. Since $(s,j,2n-b)$ is not a destroying triple, we have $\mu_{2n-b}=\lambda_j$. Since $(t,j,2n-b)$ is not a destroying triple, we have $\lambda_{2n-b}=\mu_j$. Let $x$ be an integer such that $\lambda_x=\mu_t$. If $x<2n-b$, we have a destroying triple $(x,2n-b,2n)$, a contradiction. If $x>2n-b$, choose an integer $y<a+2<2n-b, y\neq j$ such that $\mu_y \not \equiv \mu_{2n-b} \,\mathrm{mod}\, 2$. Notice that $y$ always exists. In this case, we have a destroying triple $(y,2n-b,x)$, a contradiction.
            
            \item If $\lambda_j \equiv \lambda_{2n-b}\,\mathrm{mod}\, 2$ and $\mu_j \not \equiv \mu_{2n-b} \,\mathrm{mod}\, 2$, similarly as above we have $\mu_{2n-b}=\lambda_j$. Since $\lambda_j \equiv \lambda_{2n-b} \not \equiv \lambda_{2n-i}$ for all $0 \leq i < b$ and $(t,j,2n-i)$ is not a destroying triple, we have $\lambda_{2n-i}=\mu_j$. Thus, $b=1$ and $\lambda_{2n}=\mu_j$. Let $v$ be an integer such that $\mu_v=\lambda_s$. If $v<2n-1$, we have a destroying triple $(v,2n-1,2n)$, a contradiction. If $v>2n-1$, then $v=2n$. Choose an integer $w<2n-1, w\neq s$ such that $\lambda_w \not \equiv \lambda_{2n-1} \,\mathrm{mod}\, 2$. Notice that $w$ always exists for $n\geq 3$. In this case, we have a destroying triple $(w,2n-1,2n)$ on $(\lambda,\mu)$, a contradiction.  The case when $\lambda_j \not \equiv \lambda_{2n-b} \,\mathrm{mod}\, 2$ and $\mu_j \equiv \mu_{2n-b} \,\mathrm{mod}\, 2$ can be proved in a similar way.
            
            \item If $\lambda_j \equiv \lambda_{2n-b} \,\mathrm{mod}\, 2$ and $\mu_j \equiv \mu_{2n-b} \,\mathrm{mod}\, 2$, by employing a similar argument as the previous subcase we have $b=1$, $\lambda_{2n}=\mu_j$, and $\mu_{2n}=\lambda_j$.
            If $\lambda_j=\mu_{2n-1}+1$ or $\lambda_j=\mu_{2n-1}-1$, we choose an integer $k<2n-1$ such that $\lambda_k \not \equiv \mu_{2n-1} \,\mathrm{mod}\, 2$. Notice that $k$ always exists for $n\geq 3$. In this case, we have a destroying triple $(k,2n-1,2n)$, a contradiction. If $\lambda_j \neq \mu_{2n-1}+1,\mu_{2n-1}-1$, we choose an integer $l<2n-1$ such that $\mu_l=\mu_{2n-1}-1$ or $\mu_l=\mu_{2n-1}+1$ which also satisfies $\mu_l \neq \lambda_{2n-1}$. Notice that $l$ always exists for $n\geq 3$. Then we have a destroying triple $(l,2n-1,2n)$, a contradiction.
        \end{itemize}
        
        \item Suppose that $\lambda_{a+1}\neq \mu_{a+1}$. Let $s$ be any integer such that $\lambda_s=\lambda_{a+1}-1$ or $\lambda_s=\lambda_{a+1}+1$ which also satisfies $\lambda_s \neq \mu_{a+1}$. Without loss of generality, suppose that $s>a+1$. Choose an integer $t$ such that $t<a+1$,  $\mu_t \not \equiv \mu_{a+1} \,\mathrm{mod}\, 2$ and $\mu_t\neq\lambda_{a+1},\lambda_s$. Notice that $t$ always exists since $a>1$. We have a destroying triple $(t,a+1,s)$, a contradiction. The case when $\lambda_{2n-b}\neq \mu_{2n-b}$ is similar since $1<a,b<n-1$.
        
        \item For this proof, we only consider the case when $a=n-1$. The case when $b=n-1$ is similar. Since $a=n-1$, we have $\lambda_n \not \equiv \lambda_{2n-b} \,\mathrm{mod}\, 2$. Since $b<n-1$, we have $\lambda_n \equiv \lambda_{n+1} \,\mathrm{mod}\, 2$ and $\mu_n \equiv \mu_{n+1} \,\mathrm{mod}\, 2$. By Lemma \ref{LemmaCl3} (1), we have that $\lambda_{n+1}=\mu_{n+1}$. It follows that $\lambda_n \equiv \mu_n \,\mathrm{mod}\, 2$. If $\lambda_{2n-b}=\lambda_n-1$ or $\lambda_{2n-b}=\lambda_n+1$, choose an integer $s<n$ such that $\mu_s \neq \lambda_{2n-b}$ and $\mu_s \not \equiv \mu_n \,\mathrm{mod}\, 2$. We have a destroying triple $(s,n,2n-b)$, a contradiction. Otherwise, there exists an integer $t<n$ such that $\lambda_t = \lambda_n-1$ or $\lambda_t = \lambda_n+1$ which also satisfies $\lambda_t\neq \mu_{2n-b}$. In this case, we have a destroying triple $(t,n,2n-b)$, a contradiction.
        
        \item Suppose that $\{\lambda_1,\lambda_2\} \neq \{\mu_1,\mu_2\}$. If $\lambda_2=\mu_i \in \{\mu_1,\mu_2\}$, consider a path $D \in \mathcal{D}_n$ such that $\tau_D(\mu_1)=\mu_2$. In this case, we have $\lambda(D)_2=\mathtt{u}$ and $\mu(D)_2=\mathtt{d}$. Thus, $\lambda(D) \neq \mu(D)$. The case when $\mu_2 \in \{\lambda_1,\lambda_2\}$ is similar. If $\lambda_2 \notin \{\mu_1,\mu_2\}$ and $\mu_2 \notin \{\lambda_1, \lambda_2\}$, we consider two subcases: 
        \begin{itemize}
            \item If $\mu_1=\mu_2+1$ or $\mu_1=\mu_2-1$, we choose an integer $s>2$ such that $\lambda_s \neq \mu_1,\mu_2$ and $\lambda_s \not \equiv \lambda_2 \,\mathrm{mod}\, 2$. In this case, we have a destroying triple $(1,2,s)$, a contradiction. 
            \item If $\mu_1\neq \mu_2+1,\mu_2-1$, we choose an integer $t>2$ such that $\mu_t=\mu_2+1$ or $\mu_t=\mu_2-1$ which also satisfies $\mu_t \notin \{\lambda_1,\lambda_2\}$. In this case, we have a destroying triple $(1,2,t)$, a contradiction.
        \end{itemize}
    \end{enumerate}
    This completes the proof.
\end{proof}

Now, we are ready to present the proof of Theorem \ref{ThmCl2} parts $(2)-(5)$.
\begin{proof}[Proof of Theorem \ref{ThmCl2} parts $(2)-(5)$.] \textcolor{white}{LOL}
\begin{enumerate}[label=(\arabic*)]
    \setcounter{enumi}{1}
    \item  $(\Longrightarrow)$ Let $\lambda,\mu \in S_{2n}$ with $\mathrm{par} (\lambda)=\mathrm{par} (\mu) = (n-1,n-1)$. Suppose that $\{\lambda_n, \lambda_{n+1}\}\neq\{\mu_n,\mu_{n+1}\}$. Without loss of generality, suppose that $\lambda_n \notin \{\mu_n,\mu_{n+1}\}$. If there exists a destroying triple, then we are done by Proposition \ref{Prop3}. Suppose otherwise. Let $s$ be an integer such that $\lambda_s=\lambda_n-1$ or $\lambda_s=\lambda_n+1$ with $\lambda_s\neq \mu_n, \mu_{n+1}$ and $s<n$. Such $s$ exists because $n\geq 3$ and $s\leq n+1$. Therefore, $(s,n,n+1)$ is a destroying triple, a contradiction.
    
    $(\Longleftarrow)$ Pick any $D \in \mathcal{D}_n$. Since $\mathrm{par}(\lambda)=\mathrm{par}(\mu)=(n-1,n-1)$, we have $\lambda(D)_i=\mu(D)_i=\mathtt{u}$ for all $i \in [n-1]$ and $\lambda(D)_j=\mu(D)_j=\mathtt{d}$ for all $j \in [n+2,2n]$. If $\tau_D(\lambda_n)=\lambda_{n+1}$, we have $\tau_D(\mu_n)=\mu_{n+1}$. Thus, $\lambda(D)_n=\mu(D)_n=\mathtt{u}$ and $\lambda(D)_{n+1}=\mu(D)_{n+1}=\mathtt{d}$. If $\tau_D(\lambda_n)\neq\lambda_{n+1}$, we have $\tau_D(\mu_n)\neq\mu_{n+1}$. In this case, $\lambda^{-1}\tau_D(\lambda_n)<n$ and $\lambda^{-1}\tau_D(\lambda_{n+1})>n+1$. By definition, $\lambda(D)_n=\mathtt{d}$ and $\lambda(D)_{n+1}=\mathtt{u}$. Similarly, we have $\mu(D)_n=\mathtt{d}$ and $\mu(D)_{n+1}=\mathtt{u}$. We conclude that $\lambda(D)=\mu(D)$ for all $D\in \mathcal{D}_n$, proving the claim.
    
    \item $(\Longrightarrow)$ Notice that $\lambda_{a+1}=\mu_{a+1}$ by Lemma $\ref{LemmaCl3}$ parts (2) - (3). We will prove that $\lambda_{[a]}=\mu_{[a]}$. Suppose otherwise, then there exists $p<a+1$ such that $\lambda_p\notin\mu_{[a]}$. Consider a path $D\in\mathcal{D}_n$ such that $\tau_D(\lambda_p)=\lambda_{a+1}$. In this case, we have $\lambda(D)_{a+1}=\mathtt{d}$. On the other hand, we have $\tau(\mu_{a+1})=\tau(\lambda_{a+1})=\lambda_p\notin\mu_{[a]}$, implying that $\mu(D)_{a+1}=\mathtt{u}$. It follows that $\lambda(D)\neq \mu(D)$, a contradiction. Thus, $\lambda_{[a]}=\mu_{[a]}$. By Lemma $\ref{LemmaCl3}$ part (1), we have $\lambda_i=\mu_i$ for all $i \in [a+2,2n-1-b]$. If $a<n-1$, we have $\lambda_{2n-b}=\mu_{2n-b}$ by Lemma $\ref{LemmaCl3}$ part (2). If $a=n-1$, we have $\lambda_{2n-b}=\mu_{2n-b}$ since $\lambda_{[a]}=\mu_{[a]}$. Therefore, $\lambda$ and $\mu$ are friends.
    
    \item $(\Longrightarrow)$ Without loss of generality, suppose that $a>1, b=1$. By Lemma $\ref{LemmaCl3}$ parts (2) - (3), we have $\mu_{a+1}=\lambda_{a+1}$. Similar to the previous case, we have $\lambda_{[a]}=\mu_{[a]}$. By Lemma $\ref{LemmaCl3}$ part (1), we have $\lambda_i=\mu_i$ for all $i \in [a+2, 2n-2]$. It follows that $\{\lambda_{2n-1},\lambda_{2n}\}=\{\mu_{2n-1},\mu_{2n}\}$. Obviously, there exists $\mu' \in \mathrm{fam}(\mu)$ with $\mu'_i=\mu_i$ for all $i \in [2n-2]$ and $\mu'_{2n-1}=\lambda_{2n-1}$. In this case, $\lambda$ and $\mu'$ are friends, as desired.
    
    \item $(\Longrightarrow)$ This is a direct consequence of Lemma $\ref{LemmaCl3}$ parts (1) and (4). 
\end{enumerate}
The converse of (3), (4), (5) are direct consequences of Propositions \ref{Prop1} and \ref{Prop4}.
\end{proof}

\section*{Appendix B: Table of notation and terminologies}

The following table summarises the notation and terminologies we used throughout this paper and the section and page in which they first appeared.

\begin{table}[H]
    \begin{center}
    \begin{tabular}{|c|c|c||c|c|c|}
        \hline
        Notation & Section & Page & Terminologies & Section & Page \\
        \hline
        $[a,b]$, $[n]$ & 2 & 2 & tunnel, tunneling & 2 & 3 \\
        $S_n$, $\mathcal{D}_n$ & 2 & 2 & circular representation & 2 & 3 \\
        $\sigma_k$, $\sigma_{[k]}$ & 2 & 2 & $\sigma$-paths & 3.1 & 4 \\
        $\mathtt{u}, \mathtt{d}, \mathtt{u}^k, \mathtt{d}^k$ & 2 & 2 & family & 3.2 & 5 \\
        $\tau_P$ & 2 & 3 & parity & 3.2 & 5 \\
        $h(P)$, $h_k(P)$ & 2 & 3 & friends & 3.2 & 6 \\
        $\sigma(\cdot)$ & 3.1 & 4 & circularly-connected permutation & 4 & 9 \\
        $\sim$, $\mathrm{class}(\sigma)$, $\mathrm{par}(\sigma)$  & 3.2 & 5 & block & 4 & 9 \\
        $\mathbf{U}_P$ and $\mathbf{D}_P$ & 3.3 & 8 & pairing permutation & 4.1 & 10 \\
        $\sim_Q$ & 3.3 & 8 & non-crossing & 4.1 & 10 \\
        $u_{a,k}(P)$ & 5 & 14 & unmatched and matched steps & 5 & 14 \\
        $u_{\max}^{(a)}(P)$ & 5 & 15 & destroying triple & App. A & 17 \\
        
        \hline
    \end{tabular}
    \end{center}
    \caption{Notation and terminologies in this paper.}
\end{table}

\end{document}